\input amstex
\input graphicx
\input miniltx

\documentstyle{amsppt}



\def\ep{{\varepsilon}}
\def\implies{{\ \Rightarrow\ }}
\define\lam{{\lambda}}

\define\inn#1#2{{\left\langle{#1},{#2}\right\rangle}}


\topmatter
\title {The weighted Fourier inequality, polarity, and reverse H\"older Inequality}\endtitle
\leftheadtext{}
\rightheadtext{Weighted Fourier Inequality}
\author {Ryan Berndt}\endauthor

\address {Otterbein University, 1 South Grove, Westerville, OH 43214}\endaddress
\email {rberndt\@otterbein.edu}\endemail

\subjclassyear{2010}
\subjclass
\nofrills{{\rm 2010} {\it Mathematics Subject Classification}.\usualspace}
Primary 42B10; Secondary 42A38
\endsubjclass

\keywords {Fourier transform, weights, Hausdorff-Young, reverse H\"older, polar, Mahler volume}\endkeywords


\dedicatory \enddedicatory

\abstract 
We examine the problem of the Fourier transform mapping one weighted Lebesgue space into another, by studying necessary conditions and sufficient conditions which expose an underlying geometry. In the necessary conditions, this geometry is connected to an old result of Mahler concerning the the measure of a convex body and its geometric polar being essentially reciprocal. An additional assumption, that the weights must belong to a reverse H\"older class, is used to formulate the sufficient condition. 
\endabstract

\endtopmatter

\document

In their influential and important book on weights Garc\'\i a-Cuerva and Rubio de Francia say that the following is  ``...one of the most interesting and hard problems concerning weighted inequalities." [8, p. 468]. Find necessary and sufficient conditions on weights (locally integrable, non-negative functions) $u$ and $v$ such that 

$$\left(\int_{\bold{R^n}} |\widehat f|^q\, u\right)^{1/q}\leq C\left(\int_{\bold{R^n}} |f|^p\,v\right)^{1/p}$$

\noindent for all $f\in L^1(\bold{R^n})$. 

We provide necessary conditions and sufficient conditions such that the weighted Fourier inequality holds, restricting our attention to $1<p,q<\infty$ to fix ideas and avoid special cases. The necessary condition fulfills one of our goals to expose the underlying geometry of this problem. It is known that the weighted Fourier inequality holds when an inequality involving decreasing rearrangements of $u$ and $v^{-1}$ holds [4], but by their very nature decreasing rearrangements disguise geometry. For example, the characteristic function of a square and a disk have the same decreasing rearrangement if the disk and square are equimeasurable. A single necessary and sufficient condition remains elusive, but we show that a necessary condition plus ``a little more" imply the weighted Fourier inequality.

{\it Necessity.} In order for the weighted Fourier inequality to hold, we find that the following geometric condition must be met by the weights. There is a constant $C>0$ such that for $p'=p/(p-1)$,

$$\left(\int_{E} u(x)\,dx\right)^{1/q}\left( \int_F v(x)^{-p'/p}\,dx\right)^{1/p'} \leq C,\eqno(1)$$

\noindent where $E$ is any bounded, measurable set containing the origin, and $F$ is any translation of its geometric polar. In fact, the set $E$ may be translated as well, and the sets $E$ and $F$ may even be interchanged. Further, we discuss that this class of weights is essentially the same if $E$ and $F$ are replaced with a rectanguloid and its polar, or an ellipsoid and its polar. Previously in [5], we examined this inequality where $E$ and $F$ were cubes of reciprocal measure. The current condition is a generalization because having reciprocal measure and being polar are essentially the same for cubes.

{\it Sufficiency.} We assume a weaker version of the necessary condition plus a reverse H\"older condition to prove weighted restricted weak-type estimates for the Fourier transform. This is accomplished by comparing either the measure $u$ or $v$ to Lebesgue measure and then applying the Hausdorff-Young inequality. The result actually implies restricted weak-type for a {\it range} of $p$ and $q$; we are then able to interpolate and obtain the (strong-type) weighted Fourier inequality.

\remark{Notation} If a function is nonnegative and locally integrable we call it a weight. The constant $C$ may vary at each appearance, and $p'=p/(p-1)$ is the conjugate exponent to $p$. Throughout the paper we will often write the integral of weights as measures, for example $u(E)=\int_E u$. The weight $v^{-p'/p}$ appears frequently, so we make the definition

$$w(x)=v(x)^{-p'/p}.$$

\noindent The definition of the Fourier transform we use is $\widehat f(z)=\int f(x)e^{-ix\cdot z}\,dx$ for $z\in \bold{R}^n$. All integrals are over the entire space $\bold{R}^n$ if a set is not otherwise indicated. Finally, we use the letter $C$ to indicate a constant whose value at each appearance may vary.

We assume $w$ is locally integrable for convenience. Its local integrability is actually a consequence of the local integrability of $u$ and $v$ and the weighted Fourier inequality. One can see this by first replacing $v$ with $v+\ep$ in Theorem 1 below and letting $\ep\to0$. 
\endremark

\vskip1em \noindent 1. {\bf Examples}

Before we begin the main parts of the paper, we discuss a few examples important to the problem. There are straight-forward examples of which we should be aware, as well as deeper, interesting instances of the work contained herein.

\remark{Example 1} If $u,w=v^{-p'/p}\in L^1$ then the necessary condition (1) is clearly satisfied. Also, by H\"older's inequality and the trivial estimate $|\widehat f|\leq \|f\|_1$ we have

$$\aligned
\left(\int_{\bold{R^n}} |\widehat f|^q\, u\right)^{1/q}&\leq \left(\int_{\bold{R^n}}|f|\right)\|u\|_1 \cr
&= \left(\int_{\bold{R^n}}|f|v^{1/p}v^{-1/p}\right)\|u\|_1\cr
&\leq \|u\|_1\|w\|_1\left(\int_{\bold{R^n}} |f|^p\,v\right)^{1/p}.\cr
\endaligned
$$ 
  
\noindent There are other ways to meet the necessary condition (1), but the case of integrable weights $u$ and $w$ is an important, but simple, sufficient condition.\endremark

\remark{Example 2} The classical case is the Hausdorff-Young inequality \break $\|\widehat f\|_{p'}\leq C\|f\|_p$, which holds for $p'\geq 2$. In this case $u,v\equiv 1$; the weights are bounded above and below by positive constants. We consider the connection of upper and lower bounds to the weights and the Hausdorff-Young inequality.

If $p'=q\geq 2$ and $u,v^{-1}\in L^\infty$ then there is a there a constant $C$ such that  $u(x)\leq C\leq v(x)$ for almost every $x$. Applying the Hausdorff-Young inequality we have, trivially, 

$$\aligned
\left(\int_{\bold{R^n}} |\widehat f|^q\, u\right)^{1/q} &\leq C\left(\int_{\bold{R^n}} |\widehat f|^q\, \right)^{1/q}\cr
&\leq C\left(\int_{\bold{R}^n} |f|^p\right)^{1/p}\cr
&\leq \left(\int_{\bold{R^n}} |f|^p\, v\right)^{1/p}.
\endaligned$$

Now suppose the weighted Fourier inequality holds and $u^{-1},v\in L^\infty$. The precise definition of the polar of a set is in Section 2, below, but for this example it suffices to know that for $Q$ a cube, centered at the origin, and $Q^\circ$ its polar, we have $|Q||Q^\circ|=c$ where $c$ is a fixed constant. Since $u^{-1},v\in L^\infty$, the inequality (1) implies there is a constant $C$ such that

$$|Q^\circ|^{1/q}|Q|^{1/p'}\leq\left(\int_{Q^\circ} u(x)\,dx\right)^{1/q}\left(\int_{Q} v(x)^{-p'/p}\,dx\right)^{1/p'} \leq C$$

\noindent for all cubes $Q$, centered at the origin. Because $|Q^\circ||Q|=c$, there is a constant $C$ such that $|Q|^{{1/p'}-{1/q}}\leq C$ for all $Q$. This can only happen if $p'=q$, as in the case of the Hausdorff-Young inequality. 

Therefore, if $q\geq 2$ and $u,u^{-1},v,v^{-1}\in L^\infty$ then the weighted Fourier inequality holds if and only if $p'=q$. 

\endremark

\remark{Example 3} A deeper example is the following. If we take $u\in L^\infty$ and

$$v(x)=|x|^{-\frac{p}{q}+n(p-1)}$$

\noindent then a computation shows that $u$ and $v$ satisfy the conditions of Theorem 8 so that there is a constant $C$ such that 

$$\left(\int_{\bold{R^n}} |\widehat f|^q\, u\right)^{1/q}\leq C\left(\int_{\bold{R^n}} |f|^p\,v\right)^{1/p}$$

\noindent for $p$ and $q$ satisfying: $1<p\leq q<\infty$, $p'\leq q$, $q\ne 2$.
\endremark

\vskip1em \noindent{\bf 2. Necessity}\vskip.5em

\definition{Definition} Let $E$ be a bounded set containing the origin. The {\it polar} of the set $E$ is 

$$E^\circ =\{z\in \bold{R}^n:|x\cdot z|\leq 1 \hbox{\ for all\ } x\in E\}.$$

\enddefinition

The polar of a set is a convex, centrally symmetric set ({\it i.e.} $x\in E^\circ$ implies $-x\in E^\circ$). If $E$ has non-empty interior in addition to being bounded, then the polar of $E$ is a convex body. A convex body is a compact, convex set with a non-empty interior.

A result of Mahler's [13], see also [6] and [14], is that there exist dimensional constants $d_1$ and $d_2$ such that

$$d_1 \leq |E||E^\circ| \leq  d_2$$

\noindent for all centrally symmetric, convex bodies $E$. Finding the best possible constant $d_1$ in dimensions higher than 2 is open. We are not concerned with the best constants here.

We note that if $E$ contains the origin, is bounded, and has non-empty interior, then the upper bound $d_2$ still holds|even if $E$ is not itself convex. This simple fact can be observed by noticing that in this case $E^\circ$ is a centrally symmetric convex body, so Mahler's inequality holds for $E^\circ$ and $E^{\circ\circ}$. Thus, $|E^\circ||E^{\circ\circ}| \leq d_2$. Yet, $E\subset E^{\circ\circ}$, hence $|E^\circ||E|\leq d_2$. 

A common translation construction of sets $E$ appears throughout the paper. It takes the form
 
 $$(E+\mu)^\circ+\tau.$$
 
\noindent Here,  $-\mu$ is any element of $E$ so that $E+\mu$ contains the origin, and, $\tau$ is an arbitrary element of $\bold{R}^n$ so that  $(E+\mu)^\circ+\tau$ is any translation of the polar of $E+\mu$. In other words, we imagine $E$ is moved to the origin, its polar is calculated and then translated arbitrarily. As an example, if $E$ is a cube anywhere in $\bold{R}^3$ and $-\mu$ is the center of $E$, then $(E+\mu)^\circ+\tau$ is a tetrahedron anywhere in $\bold{R}^3$.

For convenience we make the following definition.

\definition{Definition} We say that $(u,v)\in Fourier(p,q)$ if and only if there is a constant $C$ such that

$$\left(\int |\widehat f|^q\, u\right)^{1/q}\leq C\left(\int |f|^p\,v\right)^{1/p}$$

\noindent for all $f\in L^1$.

\enddefinition

Our first observation is a consequence of the properties of the Fourier transform and its modulus. We will use it frequently throughout the paper, and simply refer to it as the Translation Property. Our second observation contains necessary conditions for the weighted Fourier inequality to hold.

\proclaim{Translation Property [5]} If the weighted Fourier inequality is satisfied then all translations of the weights also satisfy the same inequality. That is,

$$\aligned
(u(x),v(x))&\in Fourier(p,q)\iff\cr (u(x+\tau_1), v(x+\tau_2))&\in Fourier(p,q)\cr
\endaligned$$

\noindent for all $\tau_1,\tau_2\in {\bold R}^n$. The constant in the Fourier inequality is unchanged for the translated weights.
\endproclaim

\proclaim{Theorem 1}  Suppose $1<p,q<\infty$. In the statements below we have the implications (i)$\implies$ (ii), and (i)$\implies$(iii)$\iff$(iv)$\implies$(ii)

\roster

\item"(i)" The weighted Fourier inequality holds, $(u,v)\in Fourier(p,q)$.

\item"(ii)" There is a constant $C$ such that for all bounded, measurable sets $E$, all $-\mu\in E$, and all $\tau\in {\bold R}^n$,

$$\left(\int_{(E+\mu)^\circ+\tau} |\widehat f|^q\, u\right)^{1/q}\leq C\left(\int_E |f|^p\,v\right)^{1/p}$$

\noindent for all $f\in L^1$ supported in $E$. The sets $E$ and $(E+\mu)^\circ+\tau$ may be interchanged in the inequality. In this case, the inequality holds for all $f$ supported in $(E+\mu)^\circ+\tau$.

\item"(iii)" There is a constant $C$ such that for all bounded, measurable sets $E$, all $-\mu\in E$, and all $\tau\in {\bold R}^n$,

$$\left(\int_E f\right)\left(\int_{(E+\mu)^\circ+\tau} u\right)^{1/q}\leq C\left(\int_E f^p\,v\right)^{1/p}$$

\noindent for all nonnegative integrable functions $f$ supported in $E$. The sets $E$ and $(E+\mu)^\circ+\tau$ may be interchanged in the inequality. In this case, the inequality holds for all $f$ supported in $(E+\mu)^\circ+\tau$.

\item"(iv)" There is a constant $C$ such that for all bounded, measurable sets $E$, all $-\mu\in E$, and all $\tau\in {\bold R}^n$,

$$\left(\int_{(E+\mu)^\circ+\tau} u\right)^{1/q}\left( \int_E v^{-p'/p}\right)^{1/p'} \leq C.$$

\noindent The sets $E$ and $(E+\mu)^\circ+\tau$ may be interchanged in the inequality.

\endroster

\endproclaim

\demo{Proof} (i)$\implies$ (ii). We suppose $f$ is supported in a bounded, measurable set $E$, and integrate over the reduced domain $(E+\mu)^\circ+\tau$ on the left hand side of the weighted Fourier inequality. If $f$ is supported in $(E+\mu)^\circ+\tau$ we integrate over the reduced domain $E$ on the left hand side.

(i)$\implies$ (iii). We suppose $f$ is nonnegative, integrable, and supported in $E$. We start with the case $0\in E$. For $z\in E^\circ$ and $x\in E$,  $\cos(x\cdot z)\geq \cos(1)$, so we have the pointwise estimate

$$\aligned
|\widehat f(z)| &= \left|\int_E f(x)e^{-ix\cdot z}\,dx\right|\cr
&\geq \int_E f(x)\cos(x\cdot z)\,dx\cr
&\geq \cos(1)\int_E f(x)\,dx.
\endaligned
$$

\noindent Here, we used the fact that the modulus of the complex number $|\widehat f(z)|$ is greater than the absolute value of its real part $|\int_E f(x)\cos(x\cdot z)|$, and this quantity is nonnegative. Thus, for nonnegative, integrable $f$ supported in sets $E$ containing the origin and $(u(z),v(x))\in Fourier(p,q)$ we have

$$\left(\int_E f(x)\,dx\right)\left(\int_{E^\circ}  u(z)\,dz\right)^{1/q}\leq C\left(\int_E f(x)^p\,v(x)\,dx\right)^{1/p}.\eqno(2)$$

\noindent For the general case,  we suppose $f$ is supported in $E$ and $-\mu\in E$ so that $0\in E+\mu$. By the Translation Property, $(u(z+\tau),v(x+\mu))\in Fourier(p,q)$ for any $\tau\in \bold{R}^n$. Thus, by (2) with weights $u(z+\tau)$ and $v(x+\mu)$ substituted for $u(z)$ and $v(x)$, we have

$$\left(\int_{E+\mu} f(x+\mu)\,dx\right)\left(\int_{(E+\mu)^\circ}  u(z+\tau)\,dz\right)^{1/q}$$
$$\leq C\left(\int_{E+\mu} |f(x+\mu)|^p\,v(x+\mu)\,dx\right)^{1/p}.$$

\noindent After a change of variables we get the inequality in (iii).

To prove the inequality where $E$ and $(E+\mu)^\circ+\tau$ are interchanged, we apply the inequality (iii), which was just established. We let $E$, $-\mu\in E$, $\tau\in \bold{R}^n$, and $f$ supported in $(E+\mu)^\circ+\tau$ be given.
We know that for any $-\mu'\in (E+\mu)^\circ+\tau$ and any $\tau'\in \bold{R}^n$ we have

$$\left(\int_{(E+\mu)^\circ+\tau} f(x)\,dx\right)\left(\int_{((E+\mu)^\circ+\tau+\mu')^\circ+\tau'}  u(z+\mu)\,dz\right)^{1/q}$$
$$\leq C\left(\int_{(E+\mu)^\circ+\tau} f(x)^p\,v(x)\,dx\right)^{1/p},$$
\noindent where we have applied the inequality in (iii) to the weight pair $(u(z+\mu),v(x))$. Now, if we choose $-\mu'$ so that $\mu'+\tau=0$ and $\tau'=0$ then the integral of the weight $u$, above, satisfies

$$\left(\int_E  u(z)\,dz\right)^{1/q}\leq \left(\int_{((E+\mu)^\circ+\tau+\mu')^\circ+\tau'}  u(z+\mu)\,dz\right)^{1/q}$$

\noindent by a change of variables and since $(E+\mu)\subset (E+\mu)^{\circ\circ}$.

(iii) $\implies$ (iv). We take $f=v^{-p'/p}\chi_E$ in the expression offered by (ii). The interchanged inequality is accomplished in a symmetric way.

(iv) $\implies$ (ii). Using the fact that $|\widehat f|\leq \|f\|_1$, H\"older's inequality, and that $f$ is supported in $E$, we write

$$\aligned
\left(\int_{(E+\mu)^\circ+\tau} |\widehat f|^q\, u\right)^{1/q}&\leq \left(\int_E f\right) u((E+\mu)^\circ+\tau))^{1/q}\cr
&\leq\left(\int_E |f|v^{1/p}v^{-1/p}\right)u((E+\mu)^\circ+\tau)^{1/q}\cr
&\leq \left(\int_E |f|^pv\right)^{1/p}w(E)^{1/p'}u((E+\mu)^\circ+\tau)^{1/q}\cr
&\leq C\left(\int_E |f|^p\,v\right)^{1/p}.\cr
\endaligned$$

\noindent Here, we applied (iv) to make the estimate $w(E)^{1/p'}u((E+\mu)^\circ+\tau)^{1/q}\leq C$. The interchanged inequality is clear and accomplished in the same way.

(iv)$\implies$ (iii) This proof is accomplished in exactly the same way as (iv) $\implies$ (ii).

\enddemo

The technique used to prove that (i) $\implies$ (iii) originates with Benedetto and Heinig [3, p. 253], where it was used in a different setting. The same technique was applied in Berndt [5] to find a necessary condition similar to the one in (iv) above; and, this technique was used very recently by De Carli, Gorbachev, and Tikhonov [7] to prove (iv) as stated.

In light of the theorem, we make the following definition of the class $N(p,q)$, to which $u$ and $v$ must necessarily belong for the weighted Fourier inequality to hold.

\definition{Definition} Suppose $1<p,q<\infty$. We say that $(u,v)\in N(p,q)$ if and only if there is a constant $C$ such that for all bounded, measurable sets $E$, all $-\mu\in E$, and all $\tau\in {\bold R}^n$

$$\left(\int_{(E+\mu)^\circ+\tau} u(x)\,dx\right)^{1/q}\left( \int_E v(x)^{-p'/p}\,dx\right)^{1/p'} \leq C\eqno(3)$$

\noindent and,

$$\left(\int_E u(x)\,dx\right)^{1/q}\left( \int_{(E+\mu)^\circ+\tau} v(x)^{-p'/p}\,dx\right)^{1/p'} \leq C\eqno(4)$$
\enddefinition

\noindent Hence, by Theorem 1, 

$$(u,v)\in Fourier(p,q) \implies (u,v)\in N(p,q).$$

We can simplify our definition of $N(p,q)$ by substituting rectanguloids $R$ and ellipsoids $S$ for arbitrary bounded, measurable sets $E$. The sets $R$ and $S$ are not assumed to have any particular orientation with respect to the coordinate axes, but are assumed to have non-empty interiors. Likewise, in order to make the comparison with these solids we only consider measurable sets $E$ with non-empty interiors.

\definition{Definition} We let $N_{\Cal E}(p,q)$ indicate the collection of all $(u,v)$ such that (3) and (4) of the definition of $N(p,q)$ hold for all bounded, measurable sets $E$ with non-empty interiors.
\enddefinition

\definition{Definition} We let $N_{\Cal R}(p,q)$ and $N_{\Cal S}(p,q)$ indicate the collection of all $(u,v)$ such that (3) and (4) of the definition of $N(p,q)$ hold for all rectanguloids $R$ and all ellipsoids $S$, respectively.
\enddefinition

For simplicity, if we ignore the translations, (3) and (4) of the definition of $N(p,q)$ state the following. There exists a constant $C$ such that

$$u(E^\circ)^{1/q}w(E)^{1/p'}\leq C\text{\ and\ } u(E)^{1/q}w(E^\circ)^{1/p'}\leq C\eqno(5)$$

\noindent for all bounded, measurable sets $E$ containing the origin. Another way to describe $N(p,q)$ is for (5) to hold for all translated weights $u(z+\tau)$ and $v(x+\mu)$. We then need only consider sets containing the origin when comparing $N(p,q)$, $N_{\Cal E}(p,q)$, $N_{\Cal S}(p,q)$, and $N_{\Cal R}(p,q)$.

\proclaim{Theorem 2} Suppose $1<p,q<\infty$ and $n\geq 2$. Each of the following statements implies the next. If $u(\sqrt{n}A)\approx u(A)$ and $w(\sqrt{n}A)\approx w(A)$ for all rectanguloids and ellipsoids $A$ then the statements are equivalent.

\roster
\item"(i)" There is a constant $C$ such that $u(R^\circ)^{1/q}w(nR)^{1/p'}\leq C$ and\break
$u(nR)^{1/q}w(R^\circ)^{1/p'}\leq C$ for all rectanguloids $R$ centered at the origin.
\item"(ii)" There is a constant $C$ such that $u(S^\circ)^{1/q}w(\sqrt{n} S)^{1/p'}\leq C$ and\break $u(\sqrt{n} S)^{1/q}w(S^\circ)^{1/p'}\leq C$ for all ellipsoids $S$ centered at the origin
\item"(iii)" There is a constant $C$ such that $u(E^\circ)^{1/q}w(E)^{1/p'}\leq C$ and\break $u(E)^{1/q}w(E^\circ)^{1/p'}\leq C$ for sets $E$ that are bounded, measurable, contain the origin, and have non-empty interiors.
\endroster
\endproclaim

\demo{Proof} In this proof, we will repeatedly use the following facts about polarity: $A\subset B$ implies $B^\circ \subset A^\circ$, $(kA)^\circ =k^{-1}A^\circ$ for $k\in \bold{R}$, and $A^{\circ\circ\circ}=A^\circ$. 

(i)$\implies$ (ii) Let $S$ be an ellipsoid in ${\bold R}^n$. By basic analytic geometry, there exists a rectangle $R$ such that

$$R\subset S \subset \sqrt{n}R$$

\noindent and by the properties of polarity

$$\frac{1}{\sqrt{n}}R^\circ\subset S^\circ\subset R^\circ.$$

\noindent The inclusions $S^\circ\subset R^\circ$ and $\sqrt{n}S\subset nR$ imply (ii).

(ii)$\implies$ (iii) The key component is a theorem of John [10], see also Ball [2, p. 13]. Given a convex, centrally symmetric body $F$ in ${\bold R}^n$ there is an ellipsoid $S$ such that $S\subset F\subset \sqrt{n}S$. The set $E$ is given, but $E^{\circ\circ}$ is a convex symmetric body so there exists an ellipsoid $S$ such that 

$$S\subset E^{\circ\circ} \subset \sqrt{n}S$$

\noindent and by the properties of polarity

$$\frac{1}{\sqrt{n}} S^\circ\subset E^{\circ} \subset S^\circ.$$

\noindent The inclusions $E^\circ\subset S^\circ$ and $E\subset E^{\circ\circ} \subset \sqrt{n}S$ imply (iii).

(iii)$\implies$ (i) Rectanguloids are examples of bounded sets with non-empty interiors, so under the hypotheses that $u(\sqrt{n}R)\approx u(R)$ and $w(\sqrt{n}R)\approx w(R)$, we can derive the inequalities in (i) simply. 
\enddemo

The sets $N_{\Cal R}(p,q)$ and $N_{\Cal S}(p,q)$ are, in general, supersets by Theorem 2. In the event that $u(\sqrt{n}A)\approx u(A)$ and $w(\sqrt{n}A)\approx w(A)$ for $A$ ellipsoids and rectanguloids, they are equal. We then have
\nobreak
$$N(p,q)\subset N_{\Cal E}(p,q)= N_{\Cal S}(p,q)=N_{\Cal R}(p,q).$$

\noindent{\bf 3. Sufficiency}\vskip1em

We use two ingredients for proving sufficiency. The first is the necessary condition and the second is a reverse H\" older condition. The full strength of the necessary condition will not be required; instead, we replace the arbitrary bounded sets $E$ with cubes $Q$ parallel to the the coordinate axes in the definition of $N(p,q)$ to create weaker classes $N_{\Cal Q}(p,q)$ and $N_{\Cal Q'}(p,q)$. The reverse H\"older condition will take the form of comparability of measures.

\definition{Definition} Suppose $1<p,q<\infty$. We say that $(u,v)\in N_{\Cal Q}(p,q)$ if and only if there is a constant $C$ such that for all cubes $Q$ with sides parallel to the coordinate axes, $-\mu\in Q$, and $\tau\in \bold{R}^n$,

$$u(Q)^{1/q}w((Q+\mu)^\circ+\tau)^{1/p'}\leq C$$

\noindent and

$$u((Q+\mu)^\circ+\tau)^{1/q}w(Q)^{1/p'}\leq C.$$

\noindent We say that $(u,v)\in N_{\Cal Q'}(p,q)$ if an only if there is a constant $C$ such that for all cubes $Q$ with sides parallel to the coordinate axes, $-\mu\in Q$, and $\tau\in \bold{R}^n$,

$$u(Q)^{1/q}|Q^\circ|\leq Cv((Q+\mu)^\circ+\tau)^{1/p}$$

\noindent and

$$u((Q+\mu)^\circ+\tau)^{1/q}|Q|\leq Cv(Q)^{1/p}.$$

\noindent We note that $N(p,q)\subset N_{\Cal R}(p,q)\subset N_{\Cal Q}(p,q)\subset N_{\Cal Q'}(p,q)$, where the last inclusion is a result of H\"older's inequality: $|E|\leq v(E)^{1/p}w(E)^{1/p'}$ for any measurable set $E$.
\enddefinition

\definition{Definition} A nonnegative measure $\mu$ is doubling if and only if there is a constant $C$ such that $\mu(2Q)\leq C \mu(Q)$ for cubes $Q$ with sides parallel to the coordinate axes. Here, $2Q$ is the two-times concentric dilation of $Q$.
\enddefinition

\definition{Definition} Suppose $\mu_1$ and $\mu_2$ are nonnegative measures where $\mu_1$ or $\mu_2$ is doubling. We say that $\mu_1$ is comparable to $\mu_2$ if and only if there is a constant $C$ and a parameter $\delta>0$ such that for all measurable sets $E$ and all cubes $Q$ where $E\subset Q$,

$$\frac{\mu_1(E)}{\mu_1(Q)}\leq C\left(\frac{\mu_2(E)}{\mu_2(Q)}\right)^{\delta}.$$

\noindent We abbreviate this as $(\mu_1,\mu_2)\in C(\delta)$.
\enddefinition

If two measures are comparable and one is doubling, then the other measure is automatically doubling (see below Theorem 3). As a consequence, if a measure is comparable to Lebesgue measure then it must be doubling. In fact, it is an $A_\infty$ weight in this case [10, p. 305]. The comparability of measures is an equivalence relation (although the parameter $\delta$ is not the same), and is equivalent to one measure and the Radon-Nikodym derivative of the other satisfying a reverse H\"older inequality for a particular exponent. 

The relationship between the parameter in comparability and the exponent in the reverse H\"older inequality is in the next theorem. We follow Garc\'\i a-Cuerva and Rubio de Francia [8, 401-402] and Grafakos [10, p. 304] closely.

\proclaim{Theorem 3}  If $\mu_1$ and $\mu_2$ are nonnegative measures where at least one of them is doubling and $(\mu_1,\mu_2)\in C(\delta)$ for some $\delta>0$, then both measures are doubling. For doubling measures, each statement below implies the next.

\roster
\item"(i)" $(\mu_1,\mu_2)\in C(\delta)$ for some $\delta>0$

\item"(ii)" There is a locally integrable, nonnegative function $\sigma$ such that $d\mu_2(x)=\sigma(x)d\mu_1(x)$ with

$$\left(\frac{1}{\mu_1(Q)}\int_Q \sigma^{1+\ep}d\mu_1\right)^{1/(1+\ep)}\leq \frac{C}{\mu_1(Q)}\int_Q \sigma d\mu_1$$

\noindent for some $\ep=\ep(\delta)>0$. The function $\ep(\delta)$ is continuous and increasing.

\item"(iii)" $(\mu_1,\mu_2)\in C(\delta)$ for $\delta=\ep/(1+\ep)$. \endroster
\endproclaim

\demo{Proof} We suppose $\mu_1$ is doubling: there is a constant $\bar C$ such that $u_1(2Q)\leq \bar C\mu_1(Q)$ for all cubes $Q$; and, $(\mu_1,\mu_2)\in C(\delta)$. Suppose $\alpha \in (0,1)$ and\break $\mu_2(E)/\mu_2(Q)\leq \alpha$ for $E\subset Q$ as in the definition of comparability. We easily see that $\mu_1(E)/\mu_1(Q)\leq C\alpha^\delta$. We let $\beta=C\alpha^\delta$ and choose $\alpha$ small enough so that $\beta\in (0,1)$. We also choose $\beta<\bar C^{-1}$, a fact we will need in a moment. Thus, there exist $\alpha,\beta\in (0,1)$ such that for any subset $E$ of a cube $Q$, $\mu_2(E)/\mu_2(Q)\leq \alpha$ implies $\mu_1(E)/\mu_1(Q)\leq \beta$. This is equivalent to saying that $\mu_1(E)/\mu_1(Q) > \beta$ implies $\mu_2(E)/\mu_2(Q)>\alpha$. Replacing $E$ with $Q$ and $Q$ with $2Q$ we have $\mu_1(2Q)<\beta^{-1}\mu_1(Q)$ implies $\mu_2(2Q)<\alpha^{-1}\mu_2(Q)$. By the way we chose $\beta$ and the fact that that $\mu_1$ is doubling, $\mu_1(2Q)<\bar C\mu_1(Q)<\beta^{-1}\mu_1(Q)$, hence $\mu_2(2Q)<\alpha^{-1}\mu_2(Q)$ for all cubes $Q$. That is, $\mu_2$ is doubling.

In the event that  $(\mu_1,\mu_2)\in C(\delta)$ and $\mu_2$ is known to be doubling, we first rewrite the implication  $\mu_2(E)/\mu_2(Q)\leq \alpha$ implies $\mu_1(E)/\mu_1(Q)\leq \beta$ by substituting $Q-E$ for $E$. This implies  $1-\mu_2(E)/\mu_2(Q)\leq \alpha$ implies $1- \mu_1(E)/\mu_1(Q)\leq \beta$. Rearranging the inequalities and considering the contrapositive we have  $\mu_1(E)/\mu_1(Q)< 1-\beta$ implies $\mu_1(E)/\mu_1(Q)< 1-\alpha$. We now proceed with the argument above to conclude $\mu_1$ is doubling.

The implication (i) $\implies$ (ii) is exactly Theorem 2.11 from [8, p. 402], except that the function $\ep(\delta)$ is left ambiguous. Grafakos [10, p. 306] carefully traces it in terms of $\alpha$ and $\beta$:

$$\ep = \frac{-\frac{1}{2}\log \beta}{\log 2^n-\log \alpha}.$$

\noindent These are any $\alpha,\beta\in (0,1)$ such that $\mu_2(E)/\mu_2(Q)\leq \alpha$ implies $\mu_1(E)/\mu_1(Q)\leq \beta$ for all measurable subsets $E$ of any cube $Q$. We saw above that we may take $\alpha \in (0,1)$ to be any number such that $C\alpha^\delta\in (0,1)$. If we choose $\alpha=(2C)^{-\delta^{-1}}$, then

$$\ep(\delta) = \frac{-\frac{1}{2}\log \frac{1}{2}}{\log 2^n+\frac{\log (2C)}{\delta}},$$

\noindent which is a nonnegative, continuous, increasing function of $\delta>0$.

Finally, the last implication (ii) $\implies$ (iii) is exactly Theorem 2.9 from [8, p. 401].

\enddemo

In our sufficiency results, we will assume that either $u$ or $v$ is comparable to Lebesgue measure, depending on the exponents $p'$ and $q$; namely, if $p'\geq q$ then (a) $(u,|\cdot|)\in C(q/p')$, and if $p'\leq q$ then (b) $(|\cdot|,v)\in C(q'/p)$. In the case $p'=q$ condition (a) or condition (b) will suffice.

Our first step is to show that the Fourier transform is restricted weak-type $(p,q)$ with respect to the measures $u$ and $v$. That is, there is a constant $C$ such that for all $\alpha>0$,

$$u(\{z:|\widehat\chi_A(z)|> \alpha\})^{1/q}\leq \frac{C}{\alpha}v(A)^{1/p}.$$

\noindent Equivalently, 

$$\|\widehat \chi_A\|_{L^{q,\infty}(u)} : = \sup_{\alpha>0} \alpha u(\{z:|\widehat\chi_A(z)|> \alpha\})^{1/q}\leq Cv(A)^{1/p}.$$

\noindent For the case where $\chi_A\not\in L^1$ we extend the definition of the Fourier transform to $L^{q,\infty}(u)$. The space $L^{q,\infty}(u)$ is a complete quasinormed space; for details about this space, see Grafakos [9, p. 2]. The method of the extension to $L^{q,\infty}(u)$ is in the proof of (iii), below.

\proclaim{Theorem 4} Let $1<p,q<\infty$ and suppose $(u,v)\in N_{\Cal Q}(p,q)$. Each of the following statements implies the next.

\roster

\item"(i)" There is a constant $C$ for all measurable sets $E$ and all cubes $Q$ such that $E\subset Q$,

$$\aligned
\frac{u(E)}{u(Q)}&\leq C\left(\frac{|E|}{|Q|}\right)^{q/p'}\hbox{\ \ if\ \ }p'\geq q\cr
\frac{|E|}{|Q|}&\leq C\left(\frac{v(E)}{v(Q)}\right)^{q'/p}\hbox{\ \ if\ \ }p'\leq q.
\endaligned$$

\item"(ii)" There is a constant $C$ for all bounded, measurable sets $E$ and $F$,

$$u(E)^{1/q}\leq C|E|^{1/s}|F|^{-1/s'}v(F)^{1/p}$$

\noindent where $s=\max\{p',q\}$.

\item"(iii)" If $\max\{p',q\}\geq 2$ the Fourier transform can be extended to functions of the form $\chi_A$, where $v(A)<\infty$, and is restricted weak-type $(p,q)$ with respect to the measures $u$ and $v$. That is, there is a constant $C$ such that 

$$u(\{z:|\widehat\chi_A(z)|> \alpha\})^{1/q}\leq \frac{C}{\alpha}v(A)^{1/p}$$

\noindent for all $\alpha>0$ and all $A$ with $v(A)<\infty$.
\endroster
\endproclaim

\demo{Proof} (i) $\implies$ (ii). We first consider the case $p'\geq q$. Suppose $E$ and $F$ are bounded, measurable subsets of ${\bold R}^n$. We let $x$ be any point in the Lebesgue set of $w$, and choose $Q$ centered at the origin. By assumption, $(u,v)\in N_{\Cal Q}(p,q)$, and since $|Q||Q^\circ|=c$,

$$\left(\frac{1}{|Q|^{q/p'}} \int_Q u\right)^{1/q}\left(\frac{1}{|Q^\circ+x|}\int_{Q^\circ+x} w\right)^{1/p'} \leq C.$$

\noindent We will shrink $Q^\circ+x$ towards $\{x\}$. As we do this, $Q$ will eventually contain $E$, so we may apply the hypothesis $\frac{u(E)}{u(Q)}\leq C\left(\frac{|E|}{|Q|}\right)^{q/p'}$. For such $Q$, we have

$$\left(\frac{1}{|E|^{q/p'}} \int_E u\right)^{1/q} \left(\frac{1}{|Q^\circ+x|}\int_{Q^\circ+x} w\right)^{1/p'} \leq C.$$

\noindent Therefore, as $Q^\circ+x\to \{x\}$ the Lebesgue differentiation theorem implies

$$\left(\frac{1}{|E|^{q/p'}}\int_E u\right)^{1/q} w(x)^{1/p'}\leq C.\eqno(6)$$

\noindent Since $x$ was an arbitrary point in the Lebesgue set of $w$, we have the same inequality for almost every $x$. We now raise both sides to the $p'$ power, integrate over $F$, and then raise both sides to the ${1/p'}$ power. Moving terms we find

$$u(E)^{1/q}w(F)^{1/p'} \leq C|E|^{1/p'}|F|^{1/p'}.\eqno(7)$$

\noindent Applying H\" older's inequality, $|F|\leq v(F)^{1/p}w(F)^{1/p'}$, we conclude

$$u(E)^{1/q}\leq C|E|^{1/p'}|F|^{-1/p}v(F)^{1/p}.$$

For the case $p'\leq q$, we suppose $E$ and $F$ are bounded, measurable sets and $Q$ is centered at the origin. Let $x$ be a point in the Lebesgue set of $u$. By the assumption $(u,v)\in N_{\Cal Q}(p,q)$ we have the inequality $u(Q+x)^{1/q}w(Q^\circ)^{1/p'}\leq C$. We apply H\"older's inequality and $|Q||Q^\circ|=c$ to find $u(Q+x)^{1/q}\leq Cv(Q^\circ)^{1/p}|Q|$ for all cubes $Q$ parallel to the coordinate axes. Dividing both sides by $|Q|^{1/q}$ we have

$$\left(\frac{1}{|Q+x|}\int_{Q+x} u\right)^{1/q}\leq C\left(\frac{1}{|Q^\circ|^{p/q'}}\int_{Q^\circ} v \right)^{1/p}.$$

\noindent As $Q+x\to \{x\}$, $Q^\circ$ eventually contains $F$, so by hypothesis $\frac{|E|}{|Q|}\leq C\left(\frac{v(E)}{v(Q)}\right)^{q'/p}$ and the Lebesgue differentiation theorem,

$$u(x)^{1/q}\leq C\left(\frac{1}{|F|^{p/q'}}\int_{F} v \right)^{1/p}.\eqno(8)$$

\noindent We raise both sides to the $q$ power, integrate over $E$, and raise to the $1/q$ power to conclude

$$u(E)^{1/q}\leq C|E|^{1/q}|F|^{-1/q'}v(F)^{1/p}.\eqno(9)$$

(ii) $\implies$ (iii). We let $A$ be a bounded, measurable  set and $E_\alpha=\{z\in \bold{R}^n:|\widehat\chi_A(z)|>\alpha\}$. Since $\chi_A\in L^1$, its Fourier transform has limit zero at infinity, hence $E_\alpha$ is bounded. Replacing $E$ with $E_\alpha$, and $F$ with $A$ in the hypothesis we have

$$u(E_\alpha)^{1/q}\leq C|E_\alpha|^{1/s}|A|^{-1/s'}v(A)^{1/p}.$$

\noindent The Hausdorff-Young inequality implies the restricted weak-type inequality
 
 $$\alpha |E_\alpha|^{1/s} \leq C|A|^{1/s'}$$
 
 \noindent for all $\alpha>0$ and for all $s\geq 2$. Thus,
 
 $$u(E_\alpha)^{1/q}\leq \frac{C}{\alpha}v(A)^{1/p}$$
 
\noindent for all bounded sets $A$. This is restricted weak-type $(p,q)$ inequality for the Fourier transform with respect to the measures $u$ and $v$, for bounded sets $A$.

Next, we extend to the case where $v(A)<\infty$, whether or not $A$ is bounded. The problem is that $\widehat \chi_A$ does not necessarily exist as the Fourier transform of an $L^1$ function. As a substitute, we take any sequence of bounded, measurable sets $A_j$ increasing to $A$ and define

$$\widehat\chi_A : = \lim_{j\to \infty} \widehat \chi_{A_j}.$$

\noindent The limit is in the $L^{q,\infty}(u)$ quasinorm. We show the limit exists, is well-defined, and extends the classical definition to sets $A$ such that $v(A)<\infty$. 

Existence is accomplished via the completeness of $L^{q,\infty}(u)$. Let $A_j$ be a sequence of bounded sets increasing to $A$. For example, $A_j$ could be the intersection of $A$ with a ball of radius $j$. We know that there is a constant $C$ such that for each $j$, and all $\alpha>0$, $\alpha u(\{z:|\widehat \chi_{A_j}(z)|>\alpha\})^{1/q}\leq Cv(A_j)^{1/p}$; or, put another way there is a constant $C$ such that

$$\|\widehat \chi_{A_j}\|_{L^{q,\infty}(u)}\leq Cv(A_j)^{1/p}\eqno(10)$$

\noindent for all $j$. 
To show the sequence $\widehat \chi_{A_j}$ is Cauchy in $L^{q,\infty}(u)$ we note that $\chi_{A_i}-\chi_{A_j}=\chi_{A_i-A_j}-\chi_{A_j-A_i}$, and since $A_i-A_j$ and $A_j-A_i$ are bounded,

$$\aligned
\|\widehat \chi_{A_i}-\widehat \chi_{A_j}\|_{L^{q,\infty}(u)}&\leq C\left(\|\widehat \chi_{A_i-A_j}\|_{L^{q,\infty}(u)}+\|\widehat \chi_{A_j-A_i}\|_{L^{q,\infty}(u)}\right)\cr
&\leq C(v(A_i-A_j)^{1/p}+v(A_j-A_i)^{1/p}).
\endaligned\eqno(11)
$$

\noindent Since $v(A)<\infty$, the dominated convergence theorem implies

$$\lim_{i,j\to\infty} v(A_i-A_j)=\int \lim_{i,j\to \infty} v(x)\chi_{A_i-A_j}(x)\,dx,$$

\noindent and, $\chi_{A_i-A_j}=\chi_{A_i}-\chi_{A_i}\chi_{A_j} \to \chi_A-\chi_A=0$. Hence, the right hand side of (11) approaches zero as $i,j\to \infty$. By completeness, $\widehat \chi_{A_j}$ must converge to some function $\phi\in L^{q,\infty}(u)$ in this space's quasinorm.

To show that $\phi$ is uniquely determined, we suppose $B_j$ is another sequence of bounded sets increasing to $A$. Again, there exists a function $\psi\in L^{q,\infty}(u)$ such that $\widehat \chi_{B_j}\to \psi$ in $L^{q,\infty}(u)$. Then,

$$\|\phi-\psi\|_{L^{q,\infty}(u)}\leq C(\|\psi-\widehat \chi_{A_j}\|_{L^{q,\infty}(u)}+\|\widehat \chi_{A_j}-\widehat \chi_{B_j}\|_{L^{q,\infty}(u)}+\|\widehat \chi_{B_j}-\psi\|_{L^{q,\infty}(u)}).$$

\noindent We already know the first and third terms on the right converge to zero as $j\to \infty$. The middle term also converges to zero since we may apply the dominated convergence theorem again and $\chi_{A_j-B_j}=\chi_{A_j}-\chi_{A_j}\chi_{B_j}\to \chi_A-\chi_A\to 0$. Therefore, $\phi$ and $\psi$ belong to the same equivalence class, because $\|f\|_{L^q,\infty(u)}=0$ implies $f=0$. We call our choice of representative function $\widehat \chi_A$. 

Finally, we show this definition is an extension of the classical definition in that it agrees with the classical definition when $A$ is unbounded but $|A|<\infty$. In this case the Fourier transform of $\chi_A$ exists classically, because $\chi_A\in L^1$. We call it ${\Cal F}$ for the purposes of this argument. Let $A_j$ be an increasing sequence of sets, converging to $A$. Then, $|{\Cal F}(z)-\widehat \chi_{A_j}(z)|\leq |A-A_j|\to 0$. By the argument above, there is a $\phi\in L^{q,\infty}(u)$ such that $\widehat \chi_{A_j}\to \phi$ in $L^{q,\infty}(u)$. Convergence in $L^{q,\infty}(u)$ implies convergence in measure [9, p. 6], which in turn implies there is a subsequence $\widehat\chi_{A_{j_k}}$ converging to $\phi$ $u$-almost everywhere. Thus, ${\Cal F}(z)=\phi(z)$ $u$-almost everywhere. 

Taking $j\to \infty$ in (10) and applying the dominated convergence theorem we get 

$$\|\widehat \chi_A\|_{{L^{q,\infty}}(u)}\leq Cv(A)^{1/p}$$

\noindent for all $A$ such that $v(A)<\infty$, completing the proof. The (extended) Fourier transform is restricted weak-type $(p,q)$ with respect to the measures $u$ and $v$.
\enddemo

\remark{Remark} For the case $p'\leq q$ we notice that we can replace the hypothesis $(u,v)\in N_{\Cal Q}(p,q)$ with the weaker hypothesis $(u,v)\in N_{\Cal Q'}(p,q)$.  We will use this fact in the proof of Theorem 8.
\endremark

With the restricted weak-type estimate in hand, we turn to proving the strong-type estimate. As usual, we will use interpolation, but the key will be that our hypotheses will imply restricted weak-type for a {\it range} of $p$ and $q$ making the interpolation possible. There are two hypotheses, the first is that the weights satisfy the necessary condition $N_{\Cal Q}(p,q)$ or $N_{\Cal Q'}(p,q)$, and the second is that one of the measures is comparable to Lebesgue measure. The next theorems provide the range of $p$ and $q$ for these hypotheses.

\proclaim{Theorem 5} Suppose $1<p_0,q_0<\infty$ and $(u,v)\in N_{\Cal Q}(p_0,q_0)$. If $p'_0\geq q_0$ and $(u,|\cdot|)\in C(q_0/p'_0)$  then $(u,v)\in N_{\Cal Q}(p,q)$ for all $(p,q)$ such that

$$\left(\frac{1}{p},\frac{1}{q}\right) \in \left\{(x,y)\in(0,1)^2: y=\frac{p'_0}{q_0}(1-x)\right\},$$

\noindent and if $p'_0\leq q_0$ and $(|\cdot|,v)\in C(q'_0/p_0)$ then $(u,v)\in N_{\Cal Q'}(p,q)$ for all $(p,q)$ such that

$$\left(\frac{1}{p},\frac{1}{q}\right) \in \left\{(x,y)\in(0,1)^2: y=-\frac{p_0}{q'_0}x+1\right\}.$$

\endproclaim

\demo{Proof} For this proof we write $w_p=v^{-p'/p}$ and $w_{p_0}=v^{-p_0'/p_0}$. We first consider the case $p'_0\geq q_0$. The hypotheses of Theorem 4 are satisfied so by (6), $v^{-1}\in L^\infty$. Thus, there is a constant $C$ such that $w_p(Q)\leq C|Q|$ for all cubes $Q$. By (7), $u(E)^{p'_0/q_0}w_{p_0}(F) \leq C|E|F|$ for all bounded sets $E$ and $F$. Let $F=\{x:v^{-1}(x)> \beta\}$ where $\beta>0$ is a fixed number so that $F$ has positive measure. If $F$ is not bounded we replace it with $F\cap B_R$, for some ball of radius $R$. Now we have, for any cube $Q$, $-\mu\in Q$, and $\tau\in {\bold R}^n$

$$\aligned
u((Q+\mu)^\circ+\tau)^{p'/q}w_p(Q)&\leq u((Q+\mu)^\circ+\tau)^{p'/q}|Q|\cr
&=u((Q+\mu)^\circ+\tau)^{p_0'/q_0}\frac{w_{p_0}(F)}{w_{p_0}(F)}|Q|\cr
&\leq |Q^\circ||Q|\frac{|F|}{w_{p_0}(F)}\cr
&\leq C
\endaligned$$

\noindent because $p'_0/q_0=p'/q$ and $|F|/w_{p_0}(F) \leq \beta^{-p'_0/p_0}$. As (7) applies to any bounded set and Lebesgue measure is translation invariant, we may interchange $Q$ and $(Q+\mu)^\circ+\tau$ in the inequality. Thus, $(u,v)\in N_{\Cal Q}(p,q)$.

The other case is $p'_0\leq q_0$. The hypotheses of Theorem 4 are satisfied. By (8) $u\in L^\infty$. Let $E=\{x:u(x)>\beta\}$ where $\beta>0$ is chosen so that $E$ has positive measure. We replace $E$ with $E\cap B_R$ for some ball of radius $R$ if $E$ is not bounded. By (9) for all bounded $F$,

$$\frac{v(F)^{1/p_0}}{|F|^{1/q'_0}}\geq \left(\frac{u(E)}{|E|}\right)^{1/q_0}\geq \beta^{1/q_0}.$$

\noindent Thus, since $p/q'=p_0/q'_0$ and taking $Q$ for $F$ above we have for any cube $Q$, $-\mu\in Q$, and $\tau\in {\bold R}^n$,

$$\aligned
u((Q+\mu)^\circ+\tau)^{p/q}|Q|^p&\leq |Q^\circ|^{p/q}|Q|^p\cr
&=C|Q|^{p/q'}\cr
&=C|Q|^{p_0/q'_0}\cr
&\leq Cv(Q).
\endaligned$$

\noindent As with the other case, since (9) applies to any bounded set and Lebesgue measure is translation invariant, we may interchange $Q$ and $(Q+\mu)^\circ+\tau$. Hence $(u,v)\in N_{\Cal Q'}(p,q)$.

\enddemo

Similarly, the comparability of measures for a particular $p$ and $q$ automatically self-improves to the comparability of measures for a range of values for $p$ and $q$.

\proclaim{Theorem 6} If $\mu_1$ and $\mu_2$ are doubling and if $(\mu_1,\mu_2)\in C(\delta)$ for some $\delta>0$ then there is a $\bar \delta > \delta$ such that

$$(\mu_1,\mu_2)\in C(\lam)\text{\ for all\ }\lambda \in [0,\bar \delta).$$
\endproclaim

\demo{Proof} The hypotheses guarantee that Gehring's Lemma, which says that the set of all $\ep$ for which the reverse H\"older inequality in (ii) of Theorem 3 holds, is an interval of the form $[0,\bar \ep)$. Gehring's lemma is true in quite general settings, see Kinnunen and Shukla [12]; in the case of doubling measures the situation is simpler and follows from the argument as in [8, p. 402].

Let $\ep(\delta)$ be the function of Theorem 3, that maps $\delta$ in the comparability inequality to $\ep(\delta)$ in the reverse H\"older inequality. The preimage $\ep^{-1}([0,\bar \ep))$ is of the form $[0,\bar \delta)$ because $\ep$ is a continuous, increasing function. Now, suppose $(\mu_1,\mu_2)\in C(\delta)$ then arithmetically $(\mu_1,\mu_2)\in C(\lambda)$ for all $\lam\leq \delta$. Therefore, $\delta$ must be in the interval $[0,\bar \delta)$ and hence $(\mu_1,\mu_2)\in C(\lambda)$ for all $\lambda \in [0,\bar \delta)$ where $\bar \delta> \delta$.
\enddemo

\proclaim{Theorem 7} If $p_0'\geq q_0$ and $(u,|\cdot|)\in C(q_0/p'_0)$ then $(u,|\cdot|)\in C(q/p')$ for all $p$ and $q$ such that
 
$$\left(\frac{1}{p},\frac{1}{q}\right) \in \left\{(x,y)\in (0,1)^2:y>-\frac{1}{\bar s}x+\frac{1}{\bar s}\right\}.$$
 
\noindent for some $\bar s>q_0/p'_0$.

If $p'_0\leq q_0$ and $(|\cdot|,v)\in C(q'_0/p_0)$ then $(|\cdot|,v)\in C(q'/p)$ for all $p$ and $q$ such that 

$$\left(\frac{1}{p},\frac{1}{q}\right) \in \left\{(x,y)\in (0,1)^2:y>-\frac{1}{\bar r}x+\frac{1}{\bar r}\right\}.$$

\noindent for some $\bar r>q'_0/p_0$.
\endproclaim

\demo{Proof} By Theorem 6, we know there exist $\bar s>q_0/p_0'$ and $\bar r>q_0'/p_0$ such that

\roster
\item"(a)"  $\frac{u(E)}{u(Q)}\leq C\left(\frac{|E|}{|Q|}\right)^{q_0/p_0'} \implies \frac{u(E)}{u(Q)}\leq C\left(\frac{|E|}{|Q|}\right)^s$ for all $s\in [0,\bar s)$
\item"(b)" $\frac{|E|}{|Q|}\leq C\left(\frac{v(E)}{v(Q)}\right)^{q_0'/p_0} \implies \frac{|E|}{|Q|}\leq C\left(\frac{v(E)}{v(Q)}\right)^{r}$ for all $r\in [0,\bar r)$.
\endroster

\noindent where (a) is for the case $p_0'\geq q_0$ and (b) is for the case $p_0'\leq q_0$. Considering the first case, we want to find all $p$ and $q$ such that $q/p'<\bar s$. Making the substitutions $x=1/p$ and $y=1/q$ we get $(1-x)/y<\bar s$, or $y>-\bar s^{-1}x +\bar s^{-1}$. The other case is handled the same way.
\enddemo

\proclaim{Corollary} If $p'_0\geq q_0$ and $(u,|\cdot|)\in C(q_0/p'_0)$ then 

$$\frac{u(E)}{u(Q)}\leq C\left(\frac{|E|}{|Q|}\right)^{q/p'}$$

\noindent for all $(1/p,1/q)$ in an open disk centered at $(1/p_0,1/q_0)$. Similarly, if $p'_0\leq q_0$ and $(|\cdot|,v)\in C(q'_0/p_0)$ then  

$$\frac{|E|}{|Q|}\leq C\left(\frac{v(E)}{v(Q)}\right)^{q'/p}$$

\noindent for all $(1/p,1/q)$ in an open disk centered at $(1/p_0,1/q_0)$.
\endproclaim

\demo{Proof} We simply note the sets of points $(1/p,1/q)$ in Theorem 7 are open sets and $(1/p_0,1/q_0)$ is an element of each in their respective cases.
\enddemo

Before we show the weighted strong-type inequality, we outline the simple idea (for the case $p'\geq q$) that might get lost in the notation. First, the two ingredients, the necessary condition and the comparability of measures, imply weighted restricted weak-type for a particular $p$ and $q$:

$$ \left\{(u,v)\in N_{\Cal Q}(p,q)\atop{(u,|\cdot|)\in C(q/p')}\right\}\implies {\text{restricted }\atop{\text{weak-type\ } (p,q)}}$$

\noindent but, together, these classes apply to a range of $p$ and $q$:

$$\left\{(u,v)\in N_{\Cal Q}(p,q)\atop{(u,|\cdot|)\in C(q/p')}\right\}\implies \left\{(u,v)\in N_{\Cal Q}(\text{line})\atop{(u,|\cdot|)\in C(\text{disk})}\right\}.$$

\noindent We are then able to deduce

$$\text{restricted weak-type for indices along the intersection of a line and a disk}$$

\noindent and therefore, by interpolation, we infer the weighted strong-type inequality for points in the interior of this intersection. We now state and prove this precisely.

\proclaim{Theorem 8} Suppose $1<p\leq q<\infty$ and $(u,v)\in N_{\Cal Q}(p,q)$. If either (a) or (b) below holds then the weighted Fourier inequality holds. That is, 

\roster
\item"(a)"  if $p'\geq q$, $p\ne 2$, and if there is a constant $C$ for all measurable sets $E$ and all cubes $Q$ such that $E\subset Q$,

$$\frac{u(E)}{u(Q)}\leq C\left(\frac{|E|}{|Q|}\right)^{q/p'}; or,$$

\item"(b)" if $p'\leq q$, $q\ne 2$, and if there is a constant $C$ for all measurable sets $E$ and all cubes $Q$ such that $E\subset Q$, 
 
$$\frac{|E|}{|Q|}\leq C\left(\frac{v(E)}{v(Q)}\right)^{q'/p}$$
 \endroster

\noindent then there is a constant $C$ such that

$$\left(\int |\widehat f|^{q}\, u\right)^{1/{q}}\leq C\left(\int |f|^{p}\,v\right)^{1/{p}}$$ 
 
 \noindent for all $f\in L^1$.
 
\endproclaim

\demo{Proof} We let $\ell$ be an open line segment. If $p'\geq q$, we let 

$$\ell=\left\{(x,y)\in (0,1)^2: y=\frac{p'}{q}(1-x)\right\}$$

\noindent if $p'\leq q$ we let

$$\ell=\left\{(x,y)\in (0,1)^2: y=-\frac{p}{q'}x+1\right\}.$$

\noindent By Theorem 5, $(u,v)\in N_{\Cal Q}(p,q)$ or $(u,v)\in N_{\Cal Q'}(p,q)$ for all $p$ and $q$ such that $(1/p,1/q)\in \ell$. By the corollary to Theorem 7, there is a disk centered at $(1/p,1/q)$ such that $u$ or $v$ is comparable to Lebesgue measure. The set $\ell\cap D$ is an open line segment in $(0,1)^2$ with $(1/p,1/q)$ at its center. We let

$$T=\left\{(x,y)\in (0,1)^2: y\leq x\right\}.$$

\noindent We note $(1/p,1/q)\in T$ by assumption, and all points $(1/p,1/q)$ in $T$ satisfy $\max\{p',q\}\geq 2$, as required by part (iii) of Theorem 4. Theorem 4 and the remark following Theorem 4 imply the Fourier transform is restricted weak-type $(s,t)$ for all $(1/s,1/t)\in \ell\cap D\cap T$. To apply interpolation, we need points to the left and right of $(1/p,1/q)$, along the line segment, for which restricted weak-type holds. That is, we need $\ell\cap D\cap T$ to be an {\it open} line segment. In the case $p'\geq q$, this occurs when $p\ne 2$ and in the case when $p'\leq q$ this occurs when $q\ne 2$. As these conditions are assured in the hypotheses, and since $p\leq q$, we may conclude by off-diagonal Marcinkiewicz interpolation [9, p. 62], that the Fourier transform is strong-type $(p,q)$. That is, the weighted Fourier inequality holds.
\enddemo

\noindent{\bf 4. A Conjecture on Mahler Volume}\vskip1em

Mahler's inequality tells us that for centrally symmetric convex bodies $E$, there are positive dimensional constants $d_1$ and $d_2$ such that

$$d_1\leq |E||E^\circ|\leq d_2.$$

\noindent We would like to generalize the quantity $|E||E^\circ|$, which is called the Mahler volume, using the weights of this paper. 

By Theorem 1 of Section 2, if $(u,v)\in Fourier(p,q)$ there is a constant $c_2$ so that for all such $E$,

$$u(E)^{1/q}w(E^\circ)^{1/p'}\leq c_2.$$

\noindent Indeed, we know that the sets $E$ and $E^\circ$ may be translated arbitrarily, while maintaining the upper bound $c_2$. Hence, there is a constant $c_2$ such that for all $\tau,\mu\in \bold{R}^n$ and all centrally symmetric convex bodies containing the origin $E$,

$$u(E+\mu)^{1/q}w(E^\circ+\tau)^{1/p'}\leq c_2.$$

At this point we would like to introduce a positive lower bound $c_1$ to the quantity in question as well; but, this trivially fails when either $u$ or $w$ is an integrable function. For example, suppose $\mu=\tau=0$ and $u$ is integrable and we increase the size of $E$. We know that $|E^\circ|$ must decrease in size, by Mahler's inequality. Hence, $u(E)^{1/q}$ has a finite maximum value, while $w(E^\circ)^{1/p'}$ approaches zero. Yet, we saw in Example 1 of Section 1 that if $u,w\in L^1$ trivially implies the weighted Fourier inequality. Therefore, our conjecture on a generalized Mahler volume for weights arising from the weighted Fourier inequality takes the following form.

\proclaim{Conjecture} Suppose $1<p,q<\infty$. If $(u,v)\in Fourier(p,q)$ and $u,w\not\in L^1$ then there exist positive constants $c_1$ and $c_2$ such that

$$c_1 \leq \sup u(E+\mu)^{1/q}w((E+\mu)^\circ+\tau)^{1/p'}\leq c_2,$$

\noindent where the supremum is taken over all symmetric convex bodies $E$ containing the origin, all $-\mu\in E$, and all $\tau\in \bold{R}^n$.

\endproclaim

To be clear, the unsolved part is the existence of the positive lower bound. We believe that proving the conjecture is important to understanding the connection between the necessary and sufficient conditions discussed in this paper.

\remark{Acknowledgements} The first part of this work was completed while the author was taking a sabbatical at the University of Western Ontario, hosted by Professor Gord Sinnamon. The author is thankful for the hospitality and helpful insights. The author is also in debt to Jo\~ ao Pedro Ramos who pointed out an error in a preprint.
\endremark

\refstyle{C}

\enddocument